\documentclass[12pt,leqno]{article}
\usepackage{amsmath,amssymb,amscd,latexsym}
\usepackage[all]{xy}
\newcommand{\C}{{\mathbb{C}}}
\newcommand{\F}{{\mathbb{F}}}
\newcommand{\oF}{\overline{\F}}
\newcommand{\Q}{{\mathbb{Q}}}
\newcommand{\R}{{\mathbb{R}}}
\newcommand{\Z}{{\mathbb{Z}}}
\newcommand{\abb}{\mathrm{ab}}
\newcommand{\car}{\mathrm{char}}
\newcommand{\ddet}{\mathrm{det}}
\newcommand{\et}{\mathrm{\acute{e}t}}
\newcommand{\id}{\mathrm{id}}
\renewcommand{\mod}{\;\mathrm{mod}\;}
\newcommand{\ord}{\mathrm{ord}}
\newcommand{\rank}{\mathrm{rank}}
\newcommand{\rk}{\mathrm{rk}\,}
\newcommand{\sgn}{\mathrm{sgn}\,}
\newcommand{\spec}{\mathrm{spec}\,}
\newcommand{\Imm}{\mathrm{Im}\,}
\newcommand{\Ker}{\mathrm{Ker}\,}
\newcommand{\RRe}{\mathrm{Re}\,}
\newcommand{\ssp}{\mathrm{sp}\,}

\newcommand{\tors}{\mathrm{tors}}
\newcommand{\tr}{\mathrm{tr}}
\newcommand{\Ah}{{\mathcal A}}
\newcommand{\Fh}{{\mathcal F}}
\newcommand{\Vh}{{\mathcal V}}
\newcommand{\ea}{\mathfrak{a}}
\newcommand{\eb}{\mathfrak{b}}
\newcommand{\ec}{\mathfrak{c}}
\newcommand{\tec}{\tilde{\ec}}
\newcommand{\ef}{\mathfrak{f}}
\newcommand{\eh}{\mathfrak{h}}
\newcommand{\eo}{\mathfrak{o}}
\newcommand{\epp}{\mathfrak{p}}
\newcommand{\eX}{{\mathcal X}}
\newcommand{\oeX}{\overline{\eX}}
\newcommand{\oH}{\bar{H}}
\newcommand{\oX}{\overline{X}}
\newcommand{\silo}{\stackrel{\sim}{\longrightarrow}}
\newcommand{\tei}{\, | \,}
\newcommand{\ent}{\;\widehat{=}\;}
\newcommand{\hullet}{\raisebox{0.03cm}{$\scriptstyle \bullet$}}
\newcommand{\halb}{\frac{1}{2}}
\newcommand{\dis}{\displaystyle}
\newtheorem{theorem}{Theorem}
\newtheorem{contheorem}[theorem]{Conditional Theorem}
\newtheorem{prop}[theorem]{Proposition}
\newtheorem{remark}[theorem]{Remark}
\newtheorem{conjecture}[theorem]{Conjecture}
\newtheorem{dic}[theorem]{Dictionary}
\newtheorem{constr}[theorem]{Construction}
\newenvironment{proof}{\noindent {\bf Proof}}{\mbox{}\hspace*{\fill}$\Box$}
\begin{document}
\title{A dynamical systems analogue of Lichtenbaum's conjectures on special values of Hasse--Weil zeta functions}
\author{Christopher Deninger}
\date{}
\maketitle
\section{Introduction} \label{s:1}

In recent years Lichtenbaum has conjectured a description for the special values of Hasse--Weil zeta functions in terms of ``Weil-\'etale cohomology''. For varieties over finite fields this cohomology theory goes back to Deligne and is easy to define. In \cite{Li1} Lichtenbaum proves his conjectures for e.g. smooth projective varieties over finite fields. See also \cite{Ge} and \cite{Ra} for interesting contributions. For flat algebraic schemes over $\spec \Z$, Lichtenbaum has ideas for the definition of Weil-\'etale cohomology but has worked them out only for spectra of number rings in \cite{Li2}. One ingredient in Lichtenbaum's formulation is the determinant of an acyclic complex of vector spaces with volume forms as it also appears in the definition of Reidemeister torsion for example.

In earlier papers we studied a class of foliated dynamical systems which had some similarities with arithmetic schemes \cite{D1}, \cite{D2}. Assuming that certain zeta regularized determinants exist we now prove an analogue of Lichtenbaum's conjectures for a particularly simple class of such dynamical systems: Using results of \'Alvarez L\'opez and Kordyukov we express the leading coefficient $\zeta^*_R (0)$ of the Ruelle zeta function $\zeta_R (s)$ at $s = 0$ in terms of analytic torsion. We then apply the Cheeger--M\"uller theorem to replace the analytic torsion by the Reidemeister torsion with respect to harmonic bases. For our dynamical systems the latter can be expressed by the same recipe as the one in Lichtenbaum's conjectures.

If the flow has a cross section a topological proof for the equality of $\zeta^*_R (0)$ with Reidemeister torsion was given by Fried a long time ago in \cite{Fr1} under somewhat different assumptions. In particular for fibrations this insight goes back to Milnor \cite{Mi}. This case where topological methods suffice corresponds to varieties over finite fields in our analogy. In the general case where the foliation has dense leaves the above rather serious analytic theories seem to be needed for the proof. 

Fried has also given a proof that $\zeta^*_R (0)$ equals analytic and hence Reidemeister torsion for the geodesic flow on certain rank one symmetric spaces, \cite{Fr2}. This theme was further developed in \cite{MS}. However, in this case which is quite different dynamically from the above, it is not possible to express the Reidemeister torsion in the terms used by Lichtenbaum in his conjecture.

This note is an extended version of a letter to Lichtenbaum. It is preliminary in several respects. For one, our class of dynamical systems is too restricted. It can surely by enlarged at the cost of replacing equalities by leading term equalities in the analysis. Secondly, we did not yet prove existence for some of the zeta regularized determinants that we use. Finally, we do not look at other integer values of $s$ apart from $s = 0$. These will require the Bismut--Zhang generalization of the Cheeger--M\"uller theorem to non-unitary local coefficients.

I would like to thank Professors Masanori Katsurada and Iekata Shiokawa very much for the invitation to Keio University where this research was done and for their hospitality. I would also like to thank Jesus \'Alvarez L\'opez and Eric Leichtnam for helpful comments.

\section{Some conjectures of Lichtenbaum}

Consider a regular scheme $\eX$ which is separated and of finite type over $\spec \Z$. For any closed point $x$ of $\eX$ the residue field $\kappa (x)$ is finite and we denote its order by $Nx$. For $\RRe s > \dim \eX$ the Hasse--Weil zeta function of $\eX$ is defined by the absolutely and locally uniformly convergent Euler product
\[
\zeta_{\eX} (s) = \prod_{x \in |\eX|} (1 - Nx^{-s})^{-1} \; .
\]
Here $|\eX|$ is the set of closed points of $\eX$. For example, for $\eX = \spec \eo_K$ where $\eo_K$ is the ring of integers in a number field we obtain the Dedekind zeta function $\zeta_K (s)$. We will assume that $\zeta_{\eX} (s)$ has an analytic continuation to $s = 0$. This is known in several cases but by no means in general.

Lichtenbaum conjectures the existence of a certain ``Weil-\'etale'' cohomology theory with and without compact supports, $H^i_c (\eX , A)$ and $H^i (\eX , A)$ for topological rings $A$. See \cite{Li1}, \cite{Li2}. It should be related to the zeta function of $\eX$ as follows.

\begin{conjecture}[Lichtenbaum] \label{t21}
Let $\eX / \Z$ be as above. Then the groups $H^i_c (\eX , \Z)$ are finitely generated and vanish for $i > 2 \dim \eX + 1$. Giving $\R$ the usual topology we have
\[
H^i_c (\eX , \Z) \otimes_{\Z} \R = H^i_c (\eX, \R) \; .
\]
Moreover, there is a canonical element $\psi$ in $H^1 (\eX,  \R)$ and the following assertions hold:\\
{\bf a} The complex
\[
\ldots \xrightarrow{D} H^i_c (\eX , \R) \xrightarrow{D} H^{i+1}_c (\eX , \R) \xrightarrow{} \ldots
\]
where $Dh = \psi \cup h$, is acyclic. Note that $D^2 = 0$ because $\deg \psi = 1$.\\[0.2cm]
{\bf b} $\ord_{s=0} \zeta_{\eX} (s) = \sum_i (-1)^i i \, \rk \, H^i_c (\eX , \Z)$. \\[0.2cm]
{\bf c} For the leading coefficient $\zeta^*_{\eX} (0)$ of $\zeta_{\eX} (s)$ in the Taylor development at $s = 0$ we have the formula:
\[
\zeta^*_{\eX} (0) = \pm \prod_i |H^i_c (\eX, \Z)_{\tors} |^{(-1)^i} / \ddet (H^{\hullet}_c (\eX , \R) , D , \ef^{\hullet}) \; .
\]
Here, $\ef^i$ is a basis of $H^i_c (\eX , \Z) / \tors$.
\end{conjecture}

{\bf Explanation} For an acyclic complex of finite dimensional $\R$-vector spaces
\begin{equation} \label{eq:22}
0 \xrightarrow{} V^0 \xrightarrow{D} V^1 \xrightarrow{D} \ldots \xrightarrow{D} V^r \xrightarrow{} 0
\end{equation}
and bases $\eb^i$ of $V^i$ a determinant $\det (V^{\hullet} , D , \eb^{\hullet})$ in $\R^*_+$ is defined as follows: The exact sequence \eqref{eq:22} induces a canonical isomorphism:
\begin{equation} \label{eq:23}
\bigotimes_i (\det V^i)^{(-1)^i} \xrightarrow{\sim} \R \; .
\end{equation}
Here $L^1 = L$ and $L^{-1} = L^*$ for one-dimensional $\R$-vector spaces $L$. Every basis $\eb^i$ determines a nonzero element of $\det V^i$. Using \eqref{eq:23} these elements determine a nonzero real number whose absolute value is denoted by $\det (V^{\hullet} , D, \eb^{\hullet})$. For bases $\ea = (w_1 , \ldots , w_n)$ and $\eb = (v_1 , \ldots , v_n)$ of a finite dimensional vector space $V$ set
\[
[\eb / \ea] = \det M
\]
where $v_i = \sum_j m_{ij} w_j$. Thus we have
\[
[\ec / \ea] = [\ec / \eb] [\eb / \ea] \; .
\]
Note that $\det (V^{\hullet} , D , \eb^{\hullet})$ is unchanged if we replace the bases $\eb^i$ with bases $\ea^i$ such that $|[\eb^i / \ea^i]| = 1$ for all $i$. Thus the $\eb^i$ could be replaced by unimodularly or orthogonally equivalent bases. In particular the determinant
\[
\det (H^{\hullet}_c (\eX , \R) , D , \ef^{\hullet})
\]
does not depend on the choice of bases $\ef^i$ of $H^i_c (\eX , \Z) / \tors$. 

A more explicit way to write down $\det (V^{\hullet} , D , \eb^{\hullet})$ is the following. For all $i$ choose bases $\ec^i$ of $D (V^{i-1})$ in $V^i$ and a linearly independent set $\tilde{\ec}^{i-1}$ of vectors in $V^{i-1}$ with $D (\tilde{\ec}^{i-1}) = \ec^i$. Then $(\ec^i , \tec^i)$ is a basis of $V^i$ since \eqref{eq:22} is acyclic and we have
\begin{equation} \label{eq:24}
\det (V^{\hullet} , D , \eb^{\hullet}) = \prod_i |[\eb^i / (\ec^i , \tec^i)]|^{(-1)^i} \; .
\end{equation}
From this the following formula is obvious:

\begin{prop} \label{t25}
Let $\ea^i$ and $\eb^i$ be bases of the $V^i$ in \eqref{eq:22}. Then we have:
\[
\det (V^{\hullet} , D , \eb^{\hullet}) = \det (V^{\hullet} , D , \ea^{\hullet}) \prod_i |[\eb^i / \ea^i]|^{(-1)^i} \; .
\]
\end{prop}

For smooth projective varieties over finite fields using the Weil-\'etale topology Lichtenbaum has proved his conjecture \ref{t21} in \cite{Li1}. See also \cite{Ge} for generalizations to the singular case. The formalism also works nicely in the study of $\zeta_{\eX} (s)$ at $s = 1/2$, c.f. \cite{Ra}. If $\eX$ is the spectrum of a number field, Lichtenbaum gave a definition of ``Weil-\'etale'' cohomology groups using the Weil group of the number field. Using the formula
\[
\zeta^*_K (0) = -\frac{hR}{w}
\]
he was able to verify his conjecture except for the vanishing of cohomology in degrees greater than three, \cite{Li2}. In fact, his cohomology does not vanish in higher degrees as was recently shown by Flach and Geisser so that some modification will be necessary. 

Let us mention that in Lichtenbaum's calculation in the number ring case the element $\psi$ corresponds to the log of the id\`ele norm viewed as a homomorphism:
\begin{equation} \label{eq:26}
\psi \ent \log \| \; \| \; .
\end{equation}

\section{Gradient flows}

We begin with some arithmetic motivation for the class of dynamical systems considered in the sequel.

Consider regular schemes $\eX$ which are separated and of finite type over $\spec \Z$. If $\eX$ has characteristic $p$ then the datum of $\eX$ over $\F_p$ is equivalent to the pair $(\eX \otimes \oF_p , \varphi)$ where $\varphi = F \otimes \id_{\oF_p}$ is the $\oF_p$-linear base extension of the absolute Frobenius morphism $F$ of $\eX$ to $\eX \otimes \oF_p$. For example, the set of closed points $|\eX|$ of $\eX$ corresponds to the set of finite $\varphi$-orbits $\eo$ on the $\oF_p$-valued points $(\eX \otimes \oF_p) (\oF_p) = \eX (\oF_p)$ of $\eX \otimes \oF_p$. We have $\log N x = |\eo| \log p$ under this correspondence. Pairs $(\eX \otimes \oF_p , \varphi)$ are roughly analogous to pairs $(M, \varphi)$ where $M$ is a smooth manifold of dimension $2 \dim \eX$ and $\varphi$ is a diffeomorphism of $M$. A better analogy would be with K\"ahler manifolds and covering self maps of degree greater than one. See remark \ref{t33} below.

The $\Z$-action on $M$ via the powers of $\varphi$ can be suspended as follows to an $\R$-action on a new space. Consider the quotient
\begin{equation} \label{eq:31}
X = M \times_{p^{\Z}} \R^*_+
\end{equation}
where the subgroup $p^{\Z}$ of $\R^*_+$ acts on $M \times \R^*_+$ as follows:
\[
p^{\nu} \cdot (m , u) = (\varphi^{-\nu} (m) , p^{\nu} u) \quad \mbox{for} \; \nu \in \Z , m \in M , u \in \R \; .
\]
The group $\R$ acts smoothly on $X$ by the formula
\[
\phi^t [m,u] = [m , e^t u] \; .
\]
Here $\phi^t$ is the diffeomorphism of $X$ corresponding to $t \in \R$. Note that the closed orbits $\gamma$ of the $\R$-action on $X$ are in bijection with the finite $\varphi$-orbits $\eo$ on $M$ in such a way that the length $l (\gamma)$ of $\gamma$ satisfies the relation $l (\gamma) = |\eo| \log p$. 

Thus in the analogy of $\eX / \F_p$ with $X$, the closed points $x$ of $\eX$ correspond to the closed orbits of the $\R$-action $\phi$ on $X$ and $Nx = p^{|\eo|}$ corresponds to $e^{l (\gamma)} = p^{|\eo|}$. Moreover if $d = \dim \eX$ then $\dim X = 2d + 1$. 

The system \eqref{eq:31} has more structure. The fibres of the natural projection of $X$ to $\R^*_+ / p^{\Z}$ form a $1$-codimensional foliation $\Fh$ (in fact a fibration). The leaves (fibres) of $\Fh$ are the images of $M$ under the immersions for every $u \in \R^*_+$ sending $m$ to $[m,u]$. In particular, the leaves are transversal to the flow lines of the $\R$-action and $\phi^t$ maps leaves to leaves for every $t \in \R$.

Now the basic idea is this: If $\eX / \Z$ is flat, there is no Frobenius and hence no analogy with a discrete time dynamical system i.e. an action of $\Z$. However one obtains a reasonable analogy with a continuous time dynamical system by the following correspondence

\newpage

\begin{dic}, {\bf part 1} \label{t32_1}
 
\em 
\begin{tabular}{p{7cm}|p{7cm}}
$\eX / \Z$ $d$-dimensional regular separated scheme of finite type & triple $(X , \Fh , \phi^t)$, where $X$ is a $2d+1$-dimensional smooth manifold with a smooth $\R$-action $\phi : \R \times X \to X$ and a one-codimensional foliation $\Fh$. The leaves of $\Fh$ should be transversal to the $\R$-orbits and every diffeomorphism $\phi^t$ should map leaves to leaves. \\ \hline
closed point $x$ of $\eX$ & closed orbit $\gamma$ of the $\R$ action  \\ \hline
Norm $Nx$ of closed point $x$ & $N \gamma = \exp l (\gamma)$ for closed orbit $\gamma$ \\ \hline
Hasse--Weil zeta function \newline
$\zeta_{\eX} (s) = \prod_{x \in |\eX|} (1 - Nx^{-s})^{-1}$ & Ruelle zeta function (see section \ref{sec:4}) $\zeta_X (s) = \prod_{\gamma} (1 - N \gamma^{-s})^{\pm 1}$ \newline
(if the product makes sense) \\ \hline
$\eX \to \spec \F_p$ & triples $(X , \Fh , \phi^t)$ where $\Fh$ is given by the fibres of an $\R$-equivariant fibration $X \to \R^*_+ / p^{\Z}$
\end{tabular}
\end{dic}

\begin{remark} \label{t33}
A more accurate analogy can be motivated as follows: If $(M , \varphi)$ is a pair consisting of a manifold with a self covering $\varphi$ of degree $\deg \varphi \ge 2$ one can form the generalized solenoid $\hat{M} = \varprojlim (\ldots \xrightarrow{\varphi} M \xrightarrow{\varphi} M \xrightarrow{} \ldots)$ , c.f. \cite{CC} 11.3. On $\hat{M}$ the map $\varphi$ becomes the shift isomorphism and as before we may consider a suspension $X = \hat{M} \times_{p^{\Z}} \R^*_+$. This example suggests that more precisely schemes $\eX / \Z$ should correspond to triples $(X , \Fh, \phi^t)$ where $X$ is a $2d + 1$ dimensional generalized solenoid which is also a foliated space with a foliation $\Fh$ by K\"ahler manifolds of complex dimension $d$. It is possible to do analysis on such generalized spaces c.f. \cite{CC}. Analogs with arithmetic are studied in more detail in \cite{D2} and \cite{Le}. However, there is no Cheeger--M\"uller theorem for solenoidal manifolds yet and therefore we have to be content with the simpler analogy above.
\end{remark}

\begin{constr} \label{t34}
Let us consider a triple $(X , \Fh , \phi^t)$ as in \ref{t32_1} above, let $Y_{\phi}$ be the vector field generated by the flow $\phi^t$ and let $T \Fh$ be the tangent bundle to the foliation. Let $\omega_{\phi}$ be the one-form on $X$ defined by
\[
\omega_{\phi} \, |_{T\Fh} = 0 \quad \mbox{and} \quad \langle \omega_{\phi} , Y_{\phi} \rangle = 1 \; .
\]
One checks that $d\omega_{\phi} = 0$ and that $\omega_{\phi}$ is $\phi^t$-invariant i.e. $\phi^{t*} \omega_{\phi} = \omega_{\phi}$ for all $t \in \R$. We may view the cohomology class $\psi = [\omega_{\phi}]$ in $H^1 (X ,\R)$ defined by $\omega_{\phi}$ as a homomorphism
\[
 \psi : \pi^{\abb}_1 (X) \longrightarrow \R \; .
\]
Its image $\Lambda \subset \R$ is called the group of periods of $(X , \Fh , \phi^t)$.
\end{constr}

It is known that $\Fh$ is a fibration if and only if $\rank\, \Lambda = 1$. In this case there is an $\R$-equivariant fibration
\[
X \longrightarrow \R / \Lambda
\]
whose fibres are the leaves of $\Fh$. Incidentally, a good reference for the dynamical systems we are considering is \cite{Fa}.

\begin{remark} \label{t35}
A closed orbit $\gamma$ in $X$ defines a free homotopy class and hence a well defined element $[\gamma]$ in $\pi^{\abb}_1 (X)$. We have
\[
\psi ([\gamma]) = l (\gamma) \; .
\]
\end{remark}

In our simple analogy the role of Lichtenbaum's Weil \'etale cohomology will be played by the usual sheaf cohomology with $\Z$ or $\R$ coefficients. This shows of course that the dynamical setting is really much simpler than the arithmetic one. 

The cohomology class $\psi$ is the analogue of Lichtenbaum's $\psi$ in our dynamical setting. Note that assertions \eqref{eq:26} and \ref{t35} correspond. Namely, in the analogy with number theory, prime ideals in a number ring correspond to closed orbits $\gamma$ with $\log N\epp$ corresponding to $l (\gamma)$. So $\psi$ corresponds to $\log \|\;\|$. 

Instead of triples $(X , \Fh , \phi^t)$ we could also have looked at triples $(X ,\omega , \phi^t)$ where $\omega$ is a closed $\phi^t$-invariant one form with $\langle \omega , Y_{\phi} \rangle = 1$. Namely, $\Ker \omega$ is then an integrable $\phi^t$-invariant subbundle, hence it comes from a $1$-codimensional foliation $\Fh$ by Frobenius' integrability theorem. Clearly $\Fh$ is transversal to the flow lines of $\phi$ and $\phi$ maps leaves of $\Fh$ to leaves.

Consider also the following motivation for our foliation setting: If the leaves of $\Fh$ are Riemann surfaces varying smoothly then one may consider smooth functions $f$ on $X$ which are meromorphic on leaves and have their divisors supported on closed orbits of the flow. For compact $X$ it follows from a formula by Ghys that we have
\[
\prod_{\gamma} \; \| f \|_{\gamma} = 1 \; .
\]
In the product $\gamma$ runs over the closed orbits and 
\[
\| f \|_{\gamma} = (N\gamma)^{-\ord_{\gamma}f}
\]
where $N\gamma = e^{l(\gamma)}$ and $\ord_{\gamma} f = \ord_z (f \, |_{F_z})$. Here $z$ is any point on $\gamma$ and $F_z$ is the leaf through $z$.\\
So the foliation setting allows for a product formula where the $N \gamma$ are not all powers of the same number. If one wants an infinitely generated $\Lambda$, one must allow the flow to have fixed points ($\ent$ infinite primes). This is discussed in \cite{Ko}.

We end this motivational section listing these and more analogies:

\begin{dic}, {\bf part 2} \label{t32_2} 

\em
\begin{tabular}{p{7cm}|p{7cm}}
Number of residue characteristics $|\car (\eX)|$ of $\eX$ & rank of period group $\Lambda$ of $(X, \Fh , \phi^t)$ \\ \hline
Weil \'etale cohomology of $\eX$ & Sheaf cohomology of $X$ \\[0.2cm] \hline
Arakelov compactification $\oeX = \eX \cup \eX_{\infty}$ where $\eX_{\infty} = \eX \otimes \R$ & Triples $(\oX , \Fh , \phi^t)$ where $\oX$ is a $2d + 1$-dimensional compactification of $X$ with a flow $\phi$ and a $1$-codimensional foliation. The flow maps leaves to leaves however $\phi$ may have fixed points on $\oX$\\ \hline
$\eX (\C) / F_{\infty}$ where $F_{\infty}$ denotes complex conjugation & set of fixed points of $\phi^t$. Note that the leaf of $\Fh$ containing a fixed point is $\phi$-invariant. \\ \hline
\end{tabular}

\begin{tabular}{p{7cm}|p{7cm}}
$\spec \kappa (x) \hookrightarrow \eX$ for $x \in |\eX|$ & Embedded circle i.e. knot $\R / l (\gamma) \Z \hookrightarrow X$ corresponding to a periodic orbit $\gamma$ (map $t + l (\gamma) \Z$ to $\phi^t (x)$ for a chosen point $x$ of $\gamma$). \\ \hline
Explicit formulas of analytic number theory & transversal index theorem for $\R$-action on $X$ and Laplacian's on the leaves of $\Fh$, c.f. \cite{D3}  \\ \hline
$- \log |d_{K / \Q}|$ & Connes' Euler characteristic $\chi_{\Fh} (X , \mu)$ for Riemann surface laminations with respect to transversal measure $\mu$ defined by $\phi^t$. \\ \hline
Product formula for number fields $\prod_v \|f\|_v = 1$ & Kopei's product formula \cite{Ko}\\ \hline
\end{tabular}
\end{dic}

\section{An analogue of Lichtenbaum's conjectures} \label{sec:4}

As motivated in the previous section we look at triples $(X , \Fh , \phi^t)$ where $X$ is a smooth compact manifold of odd dimension $2d + 1$, equipped with a one-codimensional foliation $\Fh$ and $\phi^t$ is a flow mapping leaves of $\Fh$ to leaves. Moreover we assume that the flow lines meet the leaves transversally in every point. This implies that $\phi$ has no fixed points. We thus have a decomposition 
\begin{equation}
  \label{eq:41}
  TX = T \Fh \oplus T_0 X
\end{equation}
where $T_0 X$ is the rank one bundle generated by the vector field $Y_{\phi}$ attached to the flow. Let $\omega_{\phi}$ be the closed one-form on $X$ defined in construction \ref{t34} with its cohomology class $\psi := [\omega_{\phi}]$ in $H^1 (X, \R)$. 

The role of a geometric cohomology theory like $H^{\hullet}_{\et} (X \otimes_{\F_p} \oF_p , \Q_l)$ with Frobenius action is played by the (reduced) foliation cohomology:
\[
\oH^{\hullet}_{\Fh} (X) := \Ker d_{\Fh} / \overline{\Imm d_{\Fh}} \; . 
\]
Here $(\Ah^{\hullet}_{\Fh} (X) , d_{\Fh})$ with $\Ah^i_{\Fh} (X) = \Gamma (X , \Lambda^i T^* \Fh)$ is the ``de Rham complex along the leaves'', (differentials only in the leaf direction). Moreover $\overline{\Imm d_{\Fh}}$ is the closure in the smooth topology of $\Ah^{\hullet}_{\Fh} (X)$. 

The groups $\oH^{\hullet}_{\Fh} (X)$ have a smooth linear $\R$-action $\phi^*$ induced by the flow $\phi^t$. The infinitesimal generator $\theta = \lim_{t \to 0} \frac{1}{t} (\phi^{t*} - \id)$ exists on $\oH^{\hullet}_{\Fh} (X)$. It plays a similar role as the Frobenius on \'etale or crystalline cohomology.

In general, the cohomologies $\oH^{\hullet}_{\Fh} (X)$ will be infinite dimensional Fr\'echet spaces. In order to prove results one needs a relation to harmonic forms. Assume for simplicity that $X$ and $\Fh$ are oriented in a compatible way and choose a metric $g_{\Fh}$ on $T\Fh$. This gives a Hodge scalar product on $\Ah^{\hullet}_{\Fh} (X)$. We define the leafwise Laplace operator by
\[
\Delta_{\Fh} = d_{\Fh} d^*_{\Fh} + d^*_{\Fh} d_{\Fh} \quad \mbox{on} \; \Ah^{\hullet}_{\Fh} (X) \; .
\]

Then we have by a deep theorem by \'Alvarez L\'opez and Kordyukov \cite{AK1}:
\begin{equation}
  \label{eq:42}
  \oH^{\hullet}_{\Fh} (X) = \Ker \Delta_{\Fh} \; .
\end{equation}
Note that $\Delta_{\Fh}$ is not elliptic but only elliptic along the leaves of $\Fh$. Hence the standard regularity theory of elliptic operators does not apply. 

The isomorphism (\ref{eq:42}) is a consequence of the decomposition proved in \cite{AK1}
\begin{equation}
  \label{eq:43}
  \Ah^{\hullet}_{\Fh} (X) = \Ker \Delta_{\Fh} \oplus \overline{\Imm d_{\Fh} d^*_{\Fh}} \oplus \overline{\Imm d^*_{\Fh} d_{\Fh}} \; . 
\end{equation}
Consider $\Delta_{\Fh}$ as an unbounded operator on the space $\Ah^{\hullet}_{\Fh , (2)} (X)$ of $L^2$-integrable forms along the leaves. Then $\Delta_{\Fh}$ is symmetric and its closure $\overline{\Delta}_{\Fh}$ is selfadjoint. The orthogonal projection $P$ of $\Ah^{\hullet}_{\Fh, (2)} (X)$ to $\Ker \overline{\Delta}_{\Fh}$ restricts to the projection $P$ of $\Ah^{\hullet}_{\Fh} (X)$ to $\Ker \Delta_{\Fh}$ in \eqref{eq:43}. 

In our simple analogue the role of Lichtenbaum's Weil-\'etale cohomology is played by the ordinary singular cohomology with $\Z$ or $\R$-coefficients of $X$. Note that because $X$ is compact we do not have to worry about compact supports. From the arithmetic point of view we are dealing with a very simple analogue only!

Lichtenbaum's complex is replaced by
\[
(H^{\hullet} (X, \R) , D) \quad \mbox{where}\; Dh = \psi \cup h \quad \mbox{and} \; \psi = [\omega_{\phi}] \; .
\]
Now assume that the closed orbits of the flow are non-degenerate in the sense that for every closed orbit $\gamma$ the endomorphism $T_x \phi^{l (\gamma)}$ of $T_x \Fh = T_x X / \R Y_{\phi , x}$ does not have $1$ as an eigenvalue for one or equivalently any $x$ on $\gamma$. Then it should follow from the transversal index theorem of \'Alvarez L\'opez--Kordyukov \cite{AK2} that we have:
\begin{equation}
  \label{eq:44}
  \begin{array}{rcl}
    \zeta_R (s) & := & \dis \prod_{\gamma} (1 - e^{-sl (\gamma)})^{-\varepsilon_{\gamma}} \\
 & = & \dis \prod_i \ddet_{\infty} (s \cdot \id - \theta \tei \oH^i_{\Fh} (X))^{(-1)^{i+1}} \; .
  \end{array}
\end{equation}
Here $\gamma$ runs over the closed orbits, $\varepsilon_{\gamma} = \sgn \det (1 - T_x \phi^{l (\gamma)} \tei T_x \Fh)$ for any $x \in \gamma$, the Euler product converges in some right half plane and $\det_{\infty}$ is the zeta-regularized determinant. The functions $\ddet_{\infty} (s \cdot \id - \theta \tei \oH^i_{\Fh} (X))$ should be entire.

If the closed orbits of $\phi$ are degenerate one can define a Ruelle zeta function via Fuller indices \cite{Fu} and relation \eqref{eq:44} should still hold. In the present note we assume for simplicity that the action of the flow $\phi$ on $T\Fh$ is isometric with respect to $g_{\Fh}$ and we discard the condition that the closed orbits should be non-degnerate. Then $\theta$ has pure eigenvalue spectrum with finite multiplicities on $\oH^{\hullet}_{\Fh} (X) = \Ker \Delta^{\hullet}_{\Fh}$ by \cite{DS}, Theorem 2.6. We define $\zeta_R (s)$ by the formula
\[
\zeta_R (s) = \prod_i \ddet_{\infty} (s \cdot \id - \theta \tei \oH^i_{\Fh} (X))^{(-1)^{i+1}} \; ,
\]
if the individual regularized determinants exist and define entire functions. 

The following example is due to \'Alvarez L\'opez. Let $M$ be a $\Z$-covering of a closed oriented surface $N$ of genus two. The choice of a Riemannian metric on $N$ determines a Riemannian metric on $M$ with respect to which the $\Z$-action is isometric. Fix a real number $\lambda \notin \Q$ and consider the homomorphism $\alpha : \Z \to \R / \Z$ with $\alpha (1) = \lambda \mod \Z$. Then the suspension foliation of $M$ with respect to $\alpha$ is the closed manifold $X = M \times_{\Z , \alpha} \R / \Z$ with the one-codimensional foliation $\Fh$ by the images of $M \times$ point. The leaves of $\Fh$ are diffeomorphic to $M$ and dense in $X$. The product metric on $M \times \R / \Z$ induces a metric on $X$ for which the ($1$-periodic) flow $\phi^t$ defined by $\phi^t [m,x] = [m, t+x]$ is isometric. By definition $\phi^t$ is transversal to $\Fh$ and maps leaves of $\Fh$ to leaves. It follows from \cite{AH} that the reduced leafwise cohomology group $\oH^1_{\Fh} (X)$ is infinite dimensional. Note that $\omega_{\phi} = dx$ and that the group of periods $\Lambda = \Z \oplus \lambda \Z$ has rank two.


\begin{contheorem}\label{t41}
Consider a triple $(X, \Fh, \phi^t)$ as above with $\dim X = 2d + 1$ and compatible orientations of $X$ and $\Fh$ such that the flow acts isometrically with respect to a metric $g_{\Fh}$ on $T\Fh$. Assuming that all zeta-regularized determinants in the proof below exist, we have the following assertions:\\
{\bf a} The complex
\[
\ldots \xrightarrow{D} H^i (X,\R) \xrightarrow{D} H^{i+1} (X , \R) \longrightarrow \ldots
\]
where $Dh = \psi \cup h$ with $\psi = [\omega_{\phi}]$ is acyclic.\\[0.2cm]
{\bf b} \quad $\ord_{s=0} \zeta_R (s) = \sum_i (-1)^i i \, \rk H^i (X, \Z)$.\\[0.2cm]
{\bf c} For the leading coefficient in the Taylor development at $s = 0$ we have the formula:
\[
\zeta^*_R (0) = \prod_i \, |H^i (X, \Z)_{\tors} |^{(-1)^i} / \det (H^{\hullet} (X,\R) , D, \ef^{\hullet}) \; .
\]
Here, $\ef^i$ is a basis of $H^i (X, \Z) / \tors$.
\end{contheorem}

\begin{proof}
We define a metric $g$ on $X$ by $g = g_{\Fh} + g_0$ (orthogonal sum) on $TX = T\Fh \oplus T_0 X$ with $g_0$ defined by $\| Y_{\phi , x} \| = 1$ for all $x \in X$. Then by the isometry condition the usual Laplace operator decomposes
\begin{equation}
  \label{eq:45}
  \Delta = \Delta_{\Fh} \oplus \Delta_0 \quad \mbox{and} \quad \Delta_{\Fh} \Delta_0 = \Delta_0 \Delta_{\Fh} \; .
\end{equation}
Moreover the infinitesimal generator $\theta$ of the $\phi^{t*}$-action on forms commutes with $\Delta_0 , \Delta_{\Fh}$ and $\Delta$.
The decomposition (\ref{eq:41}) gives a bigrading
\begin{equation}
  \label{eq:46}
  \Lambda^n T^* X = \bigoplus_{p+q = n} \Lambda^p T^* \Fh \otimes \Lambda^q T^*_0 X \; .
\end{equation}
We have $d = d_{\Fh} + d_0$ where $d_0$ is the exterior derivative along the foliation by the flow lines and $\Delta_0 = d_0 d^*_0 + d^*_0 d_0$. By (\ref{eq:45}) the operator $\Delta$ respects the bigrading on forms induced by (\ref{eq:46}) and we have the relation \cite{DS} proof of theorem 2.6:
\begin{equation}
  \label{eq:47}
  -\theta^2 = \Delta_0 \; .
\end{equation}
One therefore has a canonical decomposition into bidegrees:
\begin{equation}
  \label{eq:48}
  \Ker \Delta^n = \omega_{\phi} \wedge (\Ker \Delta^{n-1}_{\Fh})^{\theta = 0} \oplus (\Ker \Delta^n_{\Fh})^{\theta = 0} \; .
\end{equation}
Here $\Ker \Delta_{\Fh}$ is taken on $\Ah^{\hullet}_{\Fh} (X)$ and $H^{\theta =0} = \Ker (\theta : H \to H)$. 

Because $(H^{\hullet} (X , \R) , D) \cong (\Ker \Delta^{\hullet} , \omega_{\phi}  \wedge \, \_ \, )$, assertion {\bf a} of theorem \ref{t41} follows immediately from the decomposition (\ref{eq:48}), i.e.
the complex $(H^{\hullet} (X, \R) , D)$ is acyclic.

Next we look at the decomposition:
\[
\oH^i_{\Fh} (X) \cong \Ker \Delta^i_{\Fh} = (\Ker \Delta^i_{\Fh})^{\theta = 0} \oplus \Vh^i_0
\]
where $\Vh^i_0 =$ orthogonal complement of $(\Ker \Delta^i_{\Fh})^{\theta = 0}$ in $\Ker \Delta^i_{\Fh}$.

Then we have
\begin{equation}
  \label{eq:49}
  \begin{array}{rcl}
    \zeta_R (s) & = & \dis \prod_i \ddet_{\infty} (s - \theta \tei (\Ker \Delta^i_{\Fh})^{\theta = 0})^{(-1)^{i+1}} \prod_i \ddet_{\infty} (s - \theta \tei \Vh^i_0)^{(-1)^{i+1}} \\
 & = & \dis s^{\sum_i (-1)^{i+1} \dim (\Ker \Delta^i_{\Fh})^{\theta = 0}} \prod_i \ddet_{\infty} (s - \theta \tei \Vh^i_0)^{(-1)^{i+1}} \; .
  \end{array}
\end{equation}
This gives:
\begin{eqnarray}
  \label{eq:410}
  \ord_{s = 0} \zeta_R (s) & = & \sum_i (-1)^{i+1} \dim (\Ker \Delta^i_{\Fh})^{\theta =0} \\
& = & \sum_i (-1)^i i \dim \Ker \Delta^i \quad \mbox{by (\ref{eq:48})} \nonumber \\
& = & \sum_i (-1)^i i \dim H^i (X , \R) \nonumber \; .
\end{eqnarray}
Thus we have shown part {\bf b} of theorem \ref{t41}. 

Furthermore (\ref{eq:49}) gives for the leading coefficient:
\begin{eqnarray*}
  \zeta^*_R (0) & = & \prod_i \ddet_{\infty} (-\theta \tei \Vh^i_0)^{(-1)^{i+1}} \\
& = & \prod_i \exp (\zeta'_{-\theta_i} (0))^{(-1)^i} \\
& = & \exp \sum_i (-1)^i \zeta'_{-\theta_i} (0)
\end{eqnarray*}
where $\theta_i = (\theta \; \mbox{on} \; \Ker \Delta^i_{\Fh} \; \mbox{or} \; \Vh^i_0)$ and $\zeta_T (z) = \sum_{\lambda \in \ssp T \atop \lambda \neq 0} \lambda^{-z}$ with $\arg \lambda \in (-\pi , \pi]$.

Now, using (\ref{eq:47}) i.e. $-\theta^2 = \Delta_0$, an elementary calculation shows that
\[
\zeta'_{-\theta_i} (0) = \halb \zeta'_{\Delta^i_0 \, |_{\Ker \Delta^i_{\Fh}}} (0) \; .
\]
Note also that because of (\ref{eq:45}) we have
\[
\Delta^i_0 \, |_{\Ker \Delta^i_{\Fh}} = \Delta^i \, |_{\Ker \Delta^i_{\Fh}} \; .
\]
Thus we get
\begin{equation}
  \label{eq:411}
  \zeta^*_R (0) = \exp \sum_i (-1)^i \halb \zeta'_{\Delta^i \, |_{\Vh^i_0}} (0) \; .
\end{equation}
We claim that this is equal to the inverse of the analytic torsion introduced by Ray and Singer
\[
T (X,g) = \exp \sum_i (-1)^i \frac{i}{2} \zeta'_{\Delta^i} (0) 
\]
with $\Delta^i$ on all of $\Ah^i (X)$. 

{\bf Claim:}
\begin{equation}
  \label{eq:412}
  \zeta^*_R (0) = T (X,g)^{-1} \; .
\end{equation}

\begin{proof}
  The orthogonal decompositions
\[
\Ah^i (X) = \omega_{\phi} \wedge \Ah^{i-1}_{\Fh} (X) \oplus \Ah^i_{\Fh} (X)
\]
and
\[
\Ker \Delta^i = \omega_{\phi} \wedge (\Ker \Delta^{i-1}_{\Fh})^{\theta = 0} \oplus (\Ker \Delta^i_{\Fh})^{\theta = 0}
\]
show that because $\| \omega_{\phi , x} \| = 1$ for all $x \in X$:
\[
(\Ker \Delta^i)^{\perp} = \omega_{\phi} \wedge \Vh^{i-1} \oplus \Vh^i \; ,
\]
where $\Vh^i =$ orthogonal complement of $(\Ker \Delta^i_{\Fh})^{\theta = 0}$ in $\Ah^i_{\Fh} (X)$. Hence, if $P_i$ is the orthogonal projection to $\Ker \Delta^i$ we have:
\[
\tr (e^{-t \Delta^i} - P_i) = \tr (e^{-t \Delta^{i-1} \, |_{\Vh^{i-1}}}) + \tr (e^{-t \Delta^i \, |_{\Vh^i}}) \; .
\]
The formula
\begin{equation}
  \label{eq:413}
  \zeta_{\Delta^i} (z) = \frac{1}{\Gamma (z)} \int^{\infty}_0 t^z \tr (e^{-t \Delta^i} - P_i) \frac{dt}{t} 
\end{equation}
and similarly for $\zeta_{\Delta^i \, |_{\Vh^i}} (z)$, therefore shows that:
\[
\zeta'_{\Delta^i} (0) = \zeta'_{\Delta^i \, |_{\Vh^i}} (0) + \zeta'_{\Delta^{i-1} \, |_{\Vh^{i-1}}} (0) \; .
\]
This implies that
\[
T (X,g)^{-1} = \exp \sum_i (-1)^i \halb \zeta'_{\Delta^i \, |_{\Vh^{i}}} (0) \; .
\]
By (\ref{eq:411}) it remains to show that
\begin{equation}
  \label{eq:414}
  \sum (-1)^i \zeta'_{\Delta^i \, |_{\Vh^i_0}} (0) = \sum (-1)^i \zeta'_{\Delta^i \, |_{\Vh^i}} (0) \; .
\end{equation}
The decomposition (\ref{eq:43}) gives a decomposition:
\[
\Vh^i = \Vh^i_0 \oplus \overline{\Imm d_{\Fh} d^*_{\Fh}} \oplus \overline{\Imm d^*_{\Fh} d_{\Fh}}
\]
and we have
\[
\Delta^i \, |_{\Vh^i} = \Delta^i \, |_{\Vh^i_0} \oplus ((d_{\Fh} d^*_{\Fh})_i + \Delta^i_0) \oplus ((d^*_{\Fh} d_{\Fh})_i + \Delta^i_0) \; .
\]
By formulas like (\ref{eq:413}) it suffices for (\ref{eq:414}) to prove that
\begin{equation}
  \label{eq:415}
  \tr (e^{-t ((d_{\Fh} d^*_{\Fh})_i + \Delta^i_0)}) = \tr (e^{-t (d^*_{\Fh} d_{\Fh})_{i-1} + \Delta^{i-1}_0}) \; .
\end{equation}
The idea for this is the following: We have
\begin{equation} \label{eq:416a}
d^*_{\Fh} ((d_{\Fh} d^*_{\Fh})_i + \Delta^i_0) = ((d^*_{\Fh} d_{\Fh})_{i-1} + \Delta^{i-1}_0) d^*_{\Fh}
\end{equation}
and
\[
d^*_{\Fh} : \overline{\Imm d_{\Fh} d^*_{\Fh}} \silo \Imm d^*_{\Fh} d_{\Fh} \approx \overline{\Imm d^*_{\Fh} d_{\Fh}}
\]
is an isomorphism. So the traces in (\ref{eq:415}) should be equal because traces are invariant under conjugation. The actual proof of (\ref{eq:415}) uses a modification of a method which is fundamental for the heat equation approach to the index theorem. See also \cite{RS} p. 174. It goes as follows:\\
According to \cite{AK1} we have
\[
\lim_{\tau \to \infty} e^{-\tau \Delta_{\Fh}} = \; \mbox{orthogonal projection to} \; \Ker \Delta_{\Fh} \; .
\]
It follows that
\begin{eqnarray}
\lim_{\tau \to \infty} e^{-\tau (d_{\Fh} d^*_{\Fh})_i} & = & \mbox{orthogonal projection to} \; \Ker (d_{\Fh} d^*_{\Fh})_i \nonumber \\
 & = & 0 \; . \label{eq:416b}
\end{eqnarray}
The last equation holds because we have:
\[
\Ker (d_{\Fh} d^*_{\Fh})_i = \Ker (d_{\Fh} d^*_{\Fh} \, |_{\overline{\Imm d_{\Fh} d^*_{\Fh}}}) = \Ker (\Delta^i_{\Fh} \, |_{\overline{\Imm d_{\Fh} d^*_{\Fh}}}) = 0 \; .
\]
Noting that for $t > 0$ we are dealing with trace class operators, we get
\begin{eqnarray}
\lim_{\tau \to \infty} \tr \left( e^{-t (\tau (d_{\Fh} d^*_{\Fh})_i + \Delta^i_0)} \right) & = & \tr \left( \lim_{\tau \to \infty} e^{-t \tau (d_{\Fh} d^*_{\Fh})_i} e^{-t \Delta^i_0} \right) \nonumber \\
 & \overset{\rm \eqref{eq:416b}}{=} & 0 \; . \label{eq:416c}
\end{eqnarray}
Thus we have:
\begin{eqnarray*}
\tr (e^{-t ((d_{\Fh} d^*_{\Fh})_i + \Delta^i_0)}) & \overset{\eqref{eq:416c}}{=} & - \int^{\infty}_1 \frac{d}{d\tau} \left( \tr \, e^{-t (\tau (d_{\Fh} d^*_{\Fh})_i + \Delta^i_0)} \right) \, d\tau \\
& = & t \int^{\infty}_1 \tr \left( d_{\Fh} d^*_{\Fh} e^{-t (\tau (d_{\Fh} d^*_{\Fh})_i + \Delta^i_0)} \right) \, d\tau \\
& \overset{\eqref{eq:416a}}{=} & t \int^{\infty}_1 \tr \left( d_{\Fh} e^{-t (\tau (d^*_{\Fh} d_{\Fh})_{i-1} + \Delta^{i-1}_0)} d^*_{\Fh} \right) \, d\tau  \\
& = & t \int^{\infty}_1 \tr \left( \left( d_{\Fh} e^{- \frac{t}{2} (\tau (d^*_{\Fh} d_{\Fh})_{i-1} + \Delta^{i-1}_0)} \right) \left( e^{-\frac{t}{2} (\tau (d^*_{\Fh} d_{\Fh})_{i-1} + \Delta^{i-1}_0)} d^*_{\Fh} \right) \right) \, d\tau \; .
\end{eqnarray*}
Since the operators in brackets are both trace class we can interchange them and get:
\begin{eqnarray*}
\tr \left( e^{-t ((d_{\Fh} d^*_{\Fh})_i + \Delta^i_0)} \right) & = & t \int^{\infty}_1 \tr \left( \left( e^{-\frac{t}{2} (\tau (d^*_{\Fh} d_{\Fh})_{i-1} + \Delta^{i-1}_0)} \right) \left( d^*_{\Fh} d_{\Fh} e^{-\frac{t}{2} (\tau (d^*_{\Fh} d_{\Fh})_{i-1} + \Delta^{i-1}_0)} \right) \right) \, d\tau \\
 & = & t \int^{\infty}_1 \tr \left( d^*_{\Fh} d_{\Fh} e^{-t (\tau (d^*_{\Fh} d_{\Fh})_{i-1} + \Delta^{i-1}_0)} \right) \, d\tau \\
& = & -\int^{\infty}_1 \frac{d}{d\tau} \tr \left( e^{-t (\tau (d^*_{\Fh} d_{\Fh})_{i-1} + \Delta^{i-1}_0)} \right) \, d\tau \\
& = & \tr \left( e^{-t ((d^*_{\Fh} d_{\Fh})_{i-1} + \Delta^{i-1}_0)} \right) \; .
\end{eqnarray*}
The last equation uses the following equality which is proved in the same way as formula \eqref{eq:416c}:
\[
\lim_{\tau \to 0} \tr \left( e^{-t (\tau (d^*_{\Fh} d_{\Fh})_{i-1} + \Delta^{i-1}_0)} \right) = 0 \quad \mbox{for} \; t > 0 \; .
\]
Hence the claim i.e. formula \eqref{eq:412} is proved. 
\end{proof}

Next, we use the famous Cheeger--M\"uller theorem:
\begin{equation}
  \label{eq:416}
  T (X,g) = \tau (X,g) \qquad \mbox{c.f. \cite{Ch}, \cite{M}, \cite{M2}} \; .
\end{equation}
Here $\tau (X,g) = \tau (X,\eh_{\hullet})$ is the Reidemeister torsion with respect to the volume forms on homology given by the following bases $\eh_i$ of $H_i (X, \R)$. Choose  orthonormal bases $\eh^i$ of $\Ker \Delta^i$ with respect to the Hodge scalar product and view them as bases of $H^i (X, \R)$ via $\Ker \Delta^i \silo H^i (X, \R)$. Let $\eh_i$ be the dual base to $\eh^i$. 

We have the formula
\[
\tau (X,\eb_{\hullet}) = \tau (X , \ea_{\hullet}) \prod_i |[\eb_i / \ea_i]|^{(-1)^i}
\]
for any choices of bases $\ea_i , \eb_i$ of $H_i (X, \R)$. This is immediate from the definition of Reidemeister torsion given in \cite{M2} or \cite{RS}.

Let $\ef_{\hullet} , \ef^{\hullet}$ be dual bases of $H_{\hullet} (X, \Z) / \tors$ resp. $H^{\hullet} (X, \Z) / \tors$. Then we have in particular:
\begin{equation}
  \label{eq:417}
  \tau (X,\eh_{\hullet}) = \tau (X, \ef_{\hullet}) \prod_i |[\eh_i / \ef_i]|^{(-1)^i} \; .
\end{equation}
It follows from the definition of $\tau$ in \cite{M2} or \cite{RS} that we have
\[
\tau (X, \ef_{\hullet}) = \prod_i |H_i (X, \Z)_{\tors}|^{(-1)^i} \; .
\]
The cap-isomorphism:
\[
\_ \cap [X] : H^j (X, \Z) \silo H_{2d +1 -j} (X, \Z)
\]
now implies that
\begin{equation}
  \label{eq:418}
  \tau (X, \ef_{\hullet}) = \prod_i |H^i (X, \Z)_{\tors}|^{(-1)^{i+1}} \; .
\end{equation}
Now let us look at the acyclic complex
\[
(H^{\hullet} (X, \R) , D = \psi \cup \, \_ ) \cong (\Ker \Delta^{\hullet} , \omega_{\phi} \wedge \, \_ ) \; .
\]
Using (\ref{eq:48}), we see that it is isometrically isomorphic to 
\begin{equation}
  \label{eq:421}
  (M^{\hullet -1} \oplus M^{\hullet} , D)
\end{equation}
where $M^i = (\Ker \Delta^i_{\Fh})^{\theta = 0}$ and $D(m' , m) = (m,0)$. We may choose the orthonormal basis $\eh^i$ above to be of the form $\eh^i = (\omega_{\phi} \wedge \tilde{\eh}^{i-1}  , \tilde{\eh}^i)$ where $\tilde{\eh}^i$ is an orthonormal basis of $(\Ker \Delta^i_{\Fh})^{\theta= 0}$. For this basis, i.e. $(\tilde{\eh}^{i-1} , \tilde{\eh}^i)$ in the version (\ref{eq:421}) it is trivial from the definition (\ref{eq:24}) that
\[
|\det (H^{\hullet} (X, \R) , D , \eh^{\hullet})| = 1 \; .
\]
Using proposition \ref{t25} for $\eb^{\hullet} = \eh^{\hullet}$ and $\ea^{\hullet} = \ef^{\hullet}$ we find
\begin{equation}
  \label{eq:422}
  1 = |\det (H^{\hullet} (X, \R) , D , \ef^{\hullet})| \prod_i |[\eh^i / \ef^i]|^{(-1)^i} \; .
\end{equation}
Hence we get
\begin{eqnarray*}
  \lefteqn{\prod_i |H^i (X, \Z)_{\tors}|^{(-1)^i} / \det (H^{\hullet} (X, \R) , D, \ef^{\hullet})}\\
& \overset{(\ref{eq:422})}{=} & \prod_i |H^i (X, \Z)_{\tors} |^{(-1)^i} \prod_i |[\eh^i / \ef^i]|^{(-1)^i} \\
& \overset{(\ref{eq:418})}{=} & \tau (X, \ef_{\hullet})^{-1} \prod_i |[ \eh_i / \ef_i ]|^{(-1)^{i+1}} \quad \mbox{since} \; [\eb / \ea] = [\eb^* / \ea^*]^{-1} \\
& \overset{(\ref{eq:417})}{=} & \tau (X, \eh_{\hullet})^{-1} \\
& \overset{(\ref{eq:416})}{=} & T (X, g)^{-1} \\
& \overset{(\ref{eq:412})}{=} & \zeta^*_R (0) \; .
\end{eqnarray*}
So we get part {\bf c} of theorem \ref{t41}. Here $\zeta^*_R (0)$ is always positive. It can be different for $\zeta^*_{\eX / \Z} (0)$ because in that theory there are nonzero weights. For our simple systems, $\theta$ is skew symmetric on cohomology, hence all weights are zero. In other words, the zeroes and poles of $\zeta_R (s)$ all lie on $\RRe s = 0$.
\end{proof}

Mathematisches Institut\\
Einsteinstr. 62\\
48149 M\"unster, Germany\\
c.deninger@math.uni-muenster.de
\end{document}